\documentclass[12pt,twoside,leqno]{article}
\usepackage{amsmath}
\usepackage{amssymb}
\usepackage{amsxtra}
\usepackage{amscd}
\usepackage{amsthm}
\usepackage[mathscr]{eucal}

\theoremstyle{plain}

\newtheorem{sbthm}[subsubsection]{Theorem}

\theoremstyle{definition}

\newtheorem{sbpara}[subsubsection]{}

\begin{document}

\title{Extended period domains and algebraic groups}

\author
{Kazuya Kato
\footnote{
Department of mathematics, University of Chicago, 
Chicago, Illinois, 60637, USA.},
Chikara Nakayama
\footnote{
Department of Economics, Hitotsubashi University, 
2-1 Naka, Kunitachi, Tokyo 186-8601, Japan.},
Sampei Usui
\footnote{
Graduate School of Science, Osaka University,
Toyonaka, Osaka, 560-0043, Japan.}}

\maketitle
\renewcommand{\mathbb}{\bold}

\newcommand\Cal{\mathcal}
\newcommand\define{\newcommand}
\define\gp{\mathrm{gp}}%
\define\fs{\mathrm{fs}}%
\define\an{\mathrm{an}}%
\define\mult{\mathrm{mult}}%
\define\add{\mathrm{add}}%
\define\Ker{\mathrm{Ker}\,}%
\define\Coker{\mathrm{Coker}\,}%
\define\Hom{\mathrm{Hom}\,}%
\define\Ext{\mathrm{Ext}\,}%
\define\rank{\mathrm{rank}\,}%
\define\gr{\mathrm{gr}}%
\define\cHom{\Cal{Hom}}
\define\cExt{\Cal Ext\,}%

\define\cC{\Cal C}
\define\cD{\Cal D}
\define\cO{\Cal O}
\define\cS{\Cal S}
\define\cM{\Cal M}
\define\cG{\Cal G}
\define\cH{\Cal H}
\define\cE{\Cal E}
\define\cF{\Cal F}
\define\cN{\Cal N}
\define\fF{\frak F}
\define\fg{\frak g}
\define\fh{\frak h}
\define\Dc{\check{D}}
\define\Ec{\check{E}}

\newcommand{\N}{{\mathbb{N}}}
\newcommand{\Q}{{\mathbb{Q}}}
\newcommand{\Z}{{\mathbb{Z}}}
\newcommand{\R}{{\mathbb{R}}}
\newcommand{\C}{{\mathbb{C}}}
\newcommand{\bN}{{\mathbb{N}}}
\newcommand{\bQ}{{\mathbb{Q}}}
\newcommand{\bF}{{\mathbb{F}}}
\newcommand{\bZ}{{\mathbb{Z}}}
\newcommand{\bP}{{\mathbb{P}}}
\newcommand{\bR}{{\mathbb{R}}}
\newcommand{\bC}{{\mathbb{C}}}
\newcommand{\bS}{{\bold{S}}}
\newcommand{\bbQ}{{\bar \mathbb{Q}}}
\newcommand{\ol}[1]{\overline{#1}}
\newcommand{\too}{\longrightarrow}
\newcommand{\respect}{\rightsquigarrow}
\newcommand{\compatible}{\leftrightsquigarrow}
\newcommand{\upc}[1]{\overset {\lower 0.3ex \hbox{${\;}_{\circ}$}}{#1}}
\newcommand{\Gmlog}{\bG_{m, \log}}
\newcommand{\Gm}{\bG_m}
\newcommand{\ep}{\varepsilon}
\newcommand{\Spec}{\operatorname{Spec}}
\newcommand{\val}{{\mathrm{val}}} 
\newcommand{\n}{\operatorname{naive}}
\newcommand{\bs}{\operatorname{\backslash}}
\newcommand{\Gal}{\operatorname{{Gal}}}
\newcommand{\gal}{{\rm {Gal}}({\bar \Q}/{\Q})}
\newcommand{\galp}{{\rm {Gal}}({\bar \Q}_p/{\Q}_p)}
\newcommand{\gall}{{\rm{Gal}}({\bar \Q}_\ell/\Q_\ell)}
\newcommand{\wep}{W({\bar \Q}_p/\Q_p)}
\newcommand{\wel}{W({\bar \Q}_\ell/\Q_\ell)}
\newcommand{\Ad}{{\rm{Ad}}}
\newcommand{\BS}{{\rm {BS}}}
\newcommand{\even}{\operatorname{even}}
\newcommand{\End}{{\rm {End}}}
\newcommand{\odd}{\operatorname{odd}}
\newcommand{\GL}{\operatorname{GL}}
\newcommand{\np}{\text{non-$p$}}
\newcommand{\g}{{\gamma}}
\newcommand{\G}{{\Gamma}}
\newcommand{\Lam}{{\Lambda}}
\newcommand{\La}{{\Lambda}}
\newcommand{\lam}{{\lambda}}
\newcommand{\la}{{\lambda}}
\newcommand{\uL}{{{\hat {L}}^{\rm {ur}}}}
\newcommand{\uQp}{{{\hat \Q}_p}^{\text{ur}}}
\newcommand{\sel}{\operatorname{Sel}}
\newcommand{\dt}{{\rm{Det}}}
\newcommand{\Sig}{\Sigma}
\newcommand{\fil}{{\rm{fil}}}
\newcommand{\SL}{{\rm{SL}}}
\newcommand{\spl}{{\rm{spl}}}
\newcommand{\st}{{\rm{st}}}
\newcommand{\Isom}{{\rm {Isom}}}
\newcommand{\Mor}{{\rm {Mor}}}
\newcommand{\bg}{\bar{g}}
\newcommand{\id}{{\rm {id}}}
\newcommand{\cone}{{\rm {cone}}}
\newcommand{\al}{a}
\newcommand{\ChL}{{\cal{C}}(\La)}
\newcommand{\Image}{{\rm {Image}}}
\newcommand{\toric}{{\operatorname{toric}}}
\newcommand{\torus}{{\operatorname{torus}}}
\newcommand{\Aut}{{\rm {Aut}}}
\newcommand{\Qp}{{\mathbb{Q}}_p}
\newcommand{\barQp}{{\mathbb{Q}}_p}
\newcommand{\Qpur}{{\mathbb{Q}}_p^{\rm {ur}}}
\newcommand{\Zp}{{\mathbb{Z}}_p}
\newcommand{\Zl}{{\mathbb{Z}}_l}
\newcommand{\Ql}{{\mathbb{Q}}_l}
\newcommand{\Qlur}{{\mathbb{Q}}_l^{\rm {ur}}}
\newcommand{\F}{{\mathbb{F}}}
\newcommand{\eps}{{\epsilon}}
\newcommand{\epsLa}{{\epsilon}_{\La}}
\newcommand{\epsLaVxi}{{\epsilon}_{\La}(V, \xi)}
\newcommand{\epsOLaVxi}{{\epsilon}_{0,\La}(V, \xi)}
\newcommand{\Qplin}{{\mathbb{Q}}_p(\mu_{l^{\infty}})}
\newcommand{\otimesQplin}{\otimes_{\Qp}{\mathbb{Q}}_p(\mu_{l^{\infty}})}
\newcommand{\galFl}{{\rm{Gal}}({\bar {\Bbb F}}_\ell/{\Bbb F}_\ell)}
\newcommand{\gallur}{{\rm{Gal}}({\bar \Q}_\ell/\Q_\ell^{\rm {ur}})}
\newcommand{\galFF}{{\rm {Gal}}(F_{\infty}/F)}
\newcommand{\galFv}{{\rm {Gal}}(\bar{F}_v/F_v)}
\newcommand{\galF}{{\rm {Gal}}(\bar{F}/F)}
\newcommand{\epsVxi}{{\epsilon}(V, \xi)}
\newcommand{\epsOVxi}{{\epsilon}_0(V, \xi)}
\newcommand{\plim}{\lim_
{\scriptstyle 
\longleftarrow \atop \scriptstyle n}}
\newcommand{\sig}{{\sigma}}
\newcommand{\ga}{{\gamma}}
\newcommand{\del}{{\delta}}
\newcommand{\Vss}{V^{\rm {ss}}}
\newcommand{\Bst}{B_{\rm {st}}}
\newcommand{\Dpst}{D_{\rm {pst}}}
\newcommand{\Dcrys}{D_{\rm {crys}}}
\newcommand{\DdR}{D_{\rm {dR}}}
\newcommand{\Fin}{F_{\infty}}
\newcommand{\Kla}{K_{\lambda}}
\newcommand{\Ola}{O_{\lambda}}
\newcommand{\Mla}{M_{\lambda}}
\newcommand{\Det}{{\rm{Det}}}
\newcommand{\Sym}{{\rm{Sym}}}
\newcommand{\LaSa}{{\La_{S^*}}}
\newcommand{\cX}{{\cal {X}}}
\newcommand{\MHG}{{\frak {M}}_H(G)}
\newcommand{\tauMla}{\tau(M_{\lambda})}
\newcommand{\Fvur}{{F_v^{\rm {ur}}}}
\newcommand{\Lie}{{\rm {Lie}}}
\newcommand{\LMH}{{\rm {LMH}}}
\newcommand{\cB}{{\cal {B}}}
\newcommand{\cL}{{\cal {L}}}
\newcommand{\cW}{{\cal {W}}}
\newcommand{\fq}{{\frak {q}}}
\newcommand{\cont}{{\rm {cont}}}
\newcommand{\SC}{{SC}}
\newcommand{\Om}{{\Omega}}
\newcommand{\dR}{{\rm {dR}}}
\newcommand{\crys}{{\rm {crys}}}
\newcommand{\hatSig}{{\hat{\Sigma}}}
\newcommand{\rdet}{{{\rm {det}}}}
\newcommand{\ord}{{{\rm {ord}}}}
\newcommand{\Alb}{{{\rm {Alb}}}}
\newcommand{\BdR}{{B_{\rm {dR}}}}
\newcommand{\BdRO}{{B^0_{\rm {dR}}}}
\newcommand{\Bcrys}{{B_{\rm {crys}}}}
\newcommand{\Qw}{{\mathbb{Q}}_w}
\newcommand{\barkappa}{{\bar{\kappa}}}
\newcommand{\cP}{{\Cal {P}}}
\newcommand{\cZ}{{\Cal {Z}}}
\newcommand{\oppLa}{{\Lambda^{\circ}}}
\newcommand{\bG}{{\mathbb{G}}}
\newcommand{\br}{{{\bold r}}}
\newcommand{\triv}{{\rm{triv}}}
\newcommand{\sub}{{\subset}}
\newcommand{\LD}{{D^{\star,\mild}_{\SL(2)}}}
\newcommand{\LbD}{{D^{\star}_{\SL(2)}}}
\newcommand{\dbDv}{{D^{\star}_{\SL(2),\val}}}
\newcommand{\nspl}{{{\rm nspl}}}
\newcommand{\lval}{{[\val]}}
\newcommand{\mild}{{{\rm{mild}}}}
\newcommand{\lan}{\langle}
\newcommand{\ran}{\rangle}
\newcommand{\rar}{{{\rm {rar}}}}
\newcommand{\rac}{{{\rm {rac}}}}
\newcommand{\Rep}{{\operatorname{Rep}}}

\begin{abstract} 
 For a linear algebraic group $G$ over $\Q$, we consider the period domains $D$ classifying $G$-mixed Hodge structures, and construct the extended 
period domains 
$D_{\Sig}$. 
  In particular, we give toroidal partial compactifications of 
mixed Mumford--Tate domains, mixed Shimura varieties over $\C$, and higher Albanese manifolds. 
\end{abstract}

\renewcommand{\thefootnote}{\fnsymbol{footnote}}
\footnote[0]{Primary 14C30; 
Secondary 14D07, 32G20} 

\setcounter{section}{-1}
\section{Introduction}\label{s:intro}

For a linear algebraic group $G$ over $\Q$, we consider the period domains $D$ for  $G$-mixed Hodge structures. We construct the extended 
period domains 
$D_{\Sig}$, 
the space of nilpotent orbits. 

  In the case where $D$ is the pure Mumford--Tate domain (cf.\ Green--Griffiths--Kerr's book \cite{GGK}), 
$D_{\Sigma}$ essentially coincides with the one by Kerr--Pearlstein (\cite{KP}). 

  In the general case, 
we define $D$ in Section \ref{s:D} by modifying the definition of Shimura variety over $\C$ by Deligne 
(\cite{De}). 
  In Section \ref{s:DSig}, we introduce the extended period domain $D_{\Sig}$ and state the main results \ref{t:property}, \ref{t:main}. 
  In Section \ref{s:ex}, we explain the relation with the theory for 
the usual period domains (\cite{KNU2}), and 
as examples, the Mumford--Tate domains, 
mixed Shimura varieties, and higher Albanese manifolds.

  For the $p$-adic variant of this paper, see \cite{KK}. 

  We thank Teruhisa Koshikawa for helpful discussions. 

  We omit the details of constructions and proofs in this paper, which are to be published elsewhere. 

\medskip

\section{The period domain $D$}\label{s:D}

Let $G$ be a linear algebraic group over $\Q$.
Let $G_u$ be the unipotent radical of $G$. 
  Let $\mathrm{Rep}(G)$ be the category of 
finite-dimensional linear representations of $G$ over $\Q$. 

\subsection{$G$-mixed Hodge structures}
\begin{sbpara}\label{SCR}
Let $S_{\C/\R}$ be the Weil restriction (\cite{De}) of ${\mathbb G}_m$ from $\C$ to $\R$. 
Let $w: {\mathbb G}_{m,\R}\to S_{\C/\R}$ be the canonical homomorphism.
A linear representation of $S_{\C/\R}$ over $\R$ is equivalent to a finite-dimensional $\R$-vector space $V$ endowed with a decomposition
$V_\C:=\C\otimes_\R V=\bigoplus_{p,q\in \Z} \; V_\C^{p,q}$ 
such that $\overline{V_\C^{p,q}}=V_\C^{q,p}$ for any $p,q$. 
  The corresponding decomposition is given by 
$V_\C^{p,q}=\{v\in V_\C\;|\; [z]v= z^p{\bar z}^qv\;\text{for}\;z\in \C^\times\}.$
Here $[z]$ denotes $z$ regarded as an element of $S_{\C/\R}(\R)$. 
\end{sbpara}

\begin{sbpara}
\label{D} 
  Let $h_0: S_{\C/\R}\to (G/G_u)_\R$ be a homomorphism.
  Assume that the composite ${\mathbb G}_{m,\R}\overset{w}\to S_{\C/\R}\to (G/G_u)_\R$ is $\Q$-rational and central. 
  Assume also that for one (and hence any) lifting ${\mathbb G}_{m,\R} \to G_{\R}$ of this composite, 
the adjoint action of ${\mathbb G}_{m,\R}$ on $\Lie(G_u)_\R$ is of weight $\leq -1$. 

  Then, for any $V \in \Rep(G)$, the action of ${\mathbb G}_m$  on $V$ via a 
lifting ${\mathbb G}_m \to G$ of the above ${\mathbb G}_m\to G/G_u$  defines 
a rational increasing filtration $W$ on $V$ called the weight filtration, which is independent of the lifting. 
  
  A {\it $G$-mixed Hodge structure} ({\it $G$-MHS}, for short) is an exact $\otimes$-functor from $\mathrm{Rep}(G)$ 
to the category of $\Q$-MHS keeping the underlying vector spaces with the weight filtrations. 

  We define the {\it period domain $D$ associated to $G$ and $h_0$} as the set of all isomorphism classes of $G$-MHS whose associated homomorphism $S_{\C/\R}\to (G/G_u)_\R$ is $(G/G_u)(\R)$-conjugate to $h_0$. 
\end{sbpara}

\begin{sbpara}\label{hpure}
Let $Y$ be  the set of all isomorphism classes of exact $\otimes$-functors from 
$\Rep(G)$ to the category of triples $(V, W, F),$
where $V$ is a finite-dimensional $\Q$-vector space, $W$ is an increasing filtration on $V$ (called the weight filtration), and $F$ is a decreasing filtration on $V_\C$, preserving $V$ and $W$. 

Then $G(\C)$ acts on $Y$  by changing the Hodge filtration $F$. Let $\Dc:=G(\C)D \subset Y$. Then $\Dc$ is a $G(\C)$-orbit in $Y$, 
$D$ is a
$G(\R)G_u(\C)$-orbit 
in $Y$, and $D$ is open in $\Dc$.   $\Dc$ has a structure of a complex analytic manifold as a $G(\C)$-homogeneous space, and $D$ is an
 open submanifold of $\Dc$.

\end{sbpara}

\subsection{Polarizability} 
  For a linear algebraic group $G$, let $G'$ be the commutator algebraic subgroup.
\begin{sbpara}\label{pol} 
  Let $h_0: S_{\C/\R}\to (G/G_u)_\R$ be as in \ref{D}. 
  We say that it is {\it $\R$-polarizable} if 
$\{a\in (G/G_u)'(\R)\;|\; Ca=aC\}$ is a maximal compact subgroup of $(G/G_u)'(\R)$, 
where $C$ is the image of $i\in \C^\times = S_{\C/\R}(\R)$ by $h_0$ in $(G/G_u)(\R)$. 
\end{sbpara}

\begin{sbpara}\label{pol2}
  A relationship with the usual $\bR$-polarizability is as follows (\cite{De2} 2.11). 
  Let $h_0$ be as in \ref{D}. 
  Let $H$ be a $G$-MHS such that the associated $S_{\C/\R}\to (G/G_u)_\R$ is $\R$-polarizable. 
  Let $V \in \Rep(G)$. 
  Then for each $w\in \Z$, there is an $\bR$-bilinear form on $\gr^W_w(V)_\R$ which is stable under $(G/G_u)'$ and which polarizes $\gr^W_wH(V)$. 
     
\end{sbpara}

\begin{sbpara}\label{Gamma} We will often consider a subgroup $\Gamma$ of $G(\Q)$ satisfying the following condition.

There is a faithful $V\in \Rep(G)$ and a $\Z$-lattice $L$ in $V$ such that $L$ is stable under the action of $\Gamma$.

\end{sbpara}
\begin{sbpara} Assume that $h_0: S_{\C/\R}\to (G/G_u)_\R$ is $\R$-polarizable. 
Let $\Gamma$ be a subgroup of $G(\Q)$ satisfying the condition in \ref{Gamma}.
Then the quotient space $\Gamma \bs D$ is Hausdorff. 

\end{sbpara}

\section{Space of nilpotent orbits $D_{\Sig}$}\label{s:DSig}
We define the extended period domain $D_\Sig\supset D$ as the space of nilpotent orbits,  
and state the main results.
  We fix $G$ and $h_0$ as in \ref{D}. 
  Assume that $h_0$ is $\bR$-polarizable.

\subsection{Definition of $D_{\Sig}$}

\begin{sbpara}\label{nilp1}

  A {\it nilpotent cone} is a subset $\sig$ of $\Lie(G)_\R$ 
  satisfying the following (i)--(iii).
  
  (i) $\sig=\R_{\geq 0}N_1+\dots +\R_{\geq 0}N_n$ for some $N_1,\dots,N_n\in \Lie(G)_\R$. 
  
  (ii) For any $V\in \Rep(G)$, the image of $\sig$ under the induced map $\Lie(G)_\R\to \End_\R(V)$ consists of  nilpotent operators. 
  
  (iii) $[N,N']=0$ for any $N, N'\in \sig$. 
  \end{sbpara}
  
  \begin{sbpara}\label{nilp2}
    Let $F\in \Dc$ and let $\sig$ be a nilpotent cone.
     We say that the pair $(\sig, F)$ {\it generates a nilpotent orbit} if the following (i)--(iii) are satisfied.
     
     (i)  There is a faithful $V\in \Rep(G)$ such that the action of $\sig$ on $V_{\bR}$ is admissible with respect to $W$. 

     (ii) $N F^p\subset F^{p-1}$ for any $N\in \sig$ and $p\in \Z$.
     
     (iii) Let $N_1,\dots, N_n$ be as in (i) in \ref{nilp1}. Then $\exp(\sum_{j=1}^n z_jN_j)F\in D$ if $z_j\in \C$ and $\text{Im}(z_j)\gg 0$ ($1\leq j\leq n$). 
     
  A {\it nilpotent orbit} is a pair 
$(\sig, Z)$ of a nilpotent cone and an $\exp(\sig_\C)$-orbit in $\Dc$
satisfying that for any $F\in Z$, $(\sig, F)$ generates a nilpotent orbit. Here $\sig_\C$ denotes the $\C$-linear span of $\sig$ in $\Lie(G)_\C$. 
  
  \end{sbpara}

\begin{sbpara}\label{fan}

  A {\it weak fan $\Sig$ in $\Lie(G)$} is a nonempty set of sharp rational nilpotent cones 
satisfying the conditions that it is closed under taking faces 
and that any $\sig, \sig' \in \Sig$ coincide if they have a common interior 
point and if there is an $F\in \Dc$ such that both $(\sig, F)$ and $(\sig',F)$ generate nilpotent 
orbits. 

  Let $D_{\Sig}$ be the set of all nilpotent orbits $(\sig, Z)$ such that $\sig\in \Sig$.
  Then $D$ is naturally embedded in $D_{\Sig}$.

  Let $\Gamma$ be a subgroup of $G(\Q)$ satisfying the condition in \ref{Gamma}. 
  We say that $\Sigma$ and $\Gamma$ are {\it strongly compatible} if $\Sig$ is stable under 
the adjoint action of $\Gamma$ and if any $\sig \in \Sig$ is generated by elements whose $\exp$ in $G(\R)$ belong to $\Gamma$. 
  If this is the case, $\Gamma$ naturally acts on $D_{\Sig}$. 
\end{sbpara}

\subsection{Log mixed Hodge structures}

\begin{sbpara}
  Let $S$ be an object of the category $\cB(\log)$ (\cite{KNU2} III 1.1). We denote by $\LMH (S)$ the category of log $\Q$-mixed Hodge structures over $S$ (\cite{KNU2} III 1.3). 
  
  A {\it $G$-log mixed Hodge structure}  ({\it $G$-LMH}, for short) over $S$ is an exact $\otimes$-functor from $\Rep(G)$ 
   to $\LMH (S)$. 
  
  Let $\Gamma$ be a subgroup of $G(\Q)$ satisfying the condition in \ref{Gamma}. 
  A $G$-LMH over $S$ with a $\Gamma$-level structure is a $G$-LMH $H$ over $S$ endowed with a global section of the quotient sheaf $\Gamma\bs {\cal I}$, where $\cal I$ is the following sheaf on $S^{\log}$. For an open set $U$ of $S^{\log}$, ${\cal I}(U)$ is the set of all isomorphisms $H_\Q|_U\overset{\cong}\to \text{id}$ of $\otimes$-functors from $\Rep(G)$  to the category of local systems of $\Q$-modules over $U$. 
       
     \end{sbpara}
     \begin{sbpara}\label{type}
     Let $(G, h_0)$ be as in \ref{D}, let $\Gamma$ be a subgroup of $G(\Q)$ satisfying the condition in \ref{Gamma} and let 
  $\Sig$ be a weak fan in $\Lie(G)$ which is strongly compatible with $\Gamma$. 
  A $G$-LMH $H$ over $S$ with a $\Gamma$-level structure $\lambda$  is said to be {\it of $\Sig$-type} if the following (i) and (ii) are satisfied for any $s\in S$ and any $t\in s^{\log}$.
  Take a
 $\otimes$-isomorphism $\tilde \lambda_t: H_{\Q,t}\cong  \;\text{id}$ which belongs to $\lambda_t$.

  (i) There is a $\sig\in \Sig$ such that the logarithm of the action of $\Hom((M_S/\cO^\times_S)_s, \N)\subset \pi_1(s^{\log})$ on $H_{\Q,t}$ is contained, via $\tilde \lambda_t$, in $\sigma \subset \Lie(G)_\R$.

  (ii) Let $\sig\in \Sig$ be the smallest cone satisfying (i). Let  $a: \cO_{S,t}^{\log}\to \C$ be a ring homomorphism which induces the evaluation $\cO_{S,s}\to \C$ at $s$ and consider the element $F: V\mapsto {\tilde \lambda}_t a(H(V))$ of $Y$ (\ref{hpure}). Then this element belongs to $\Dc$ and $(\sig, F)$ generates a nilpotent orbit (\ref{nilp2}).
  
  If $(H, \lambda)$ is of $\Sig$-type, we have a map $S \to \Gamma \bs D_{\Sig}$, called the {\it period map} associated to $(H, \lambda)$,  which sends $s\in S$ to the class of the nilpotent orbit $(\sig, Z)\in D_{\Sig}$ which is obtained in (ii) in \ref{type}. 
  
  \end{sbpara}

\begin{sbpara}
Let $(G,h_0, \Gamma, \Sig)$ be as in \ref{type}. We endow $\Gamma \bs D_{\Sig}$ with a topology, a sheaf of rings $\cO$ over $\C$ and a log structure $M$ defined as follows. 
The topology  is the strongest topology for which the period map $S \to \Gamma \bs D_{\Sig}$ is continuous for any $(S, H, \lambda)$, where $S$ is an object of $\cB(\log)$, $H$ is a $G$-LMH on $S$, and $\lambda$ is a $\Gamma$-level structure which is of $\Sig$-type. For an open set $U$
of $\Gamma \bs D_{\Sig}$, $\cO(U)$ (resp.\ $M(U)$) is the set of all $\C$-valued functions $f$ on $U$ such that for any $(S, H, \lambda)$ as above with the period map $\phi: S \to \Gamma \bs D_{\Sig}$, the pullback of $f$ on $U':=\phi^{-1}(U)$ belongs to the image of $\cO_{U'}$ (resp.\ $M_{U'}$) in the sheaf of $\C$-valued functions on $U'$.

These structures of $\Gamma \bs D_{\Sig}$ are defined also by defining spaces $E_{\sig}$ ($\sig\in \Sig$) in a similar way as  \cite{KNU2} III. We get the same structures when we use only $S=E_{\sig}$ for $\sig \in \Sig$ and the universal objects $(H, \lambda)$ over $E_{\sig}$ in the above definitions of the structures.
\end{sbpara}

\subsection{Main results}
Let $(G, h_0, \Gamma, \Sig)$ be as in \ref{type}. Assume that $h_0$ is $\R$-polarizable. 
\begin{sbthm}
\label{t:property}
$(1)$ $\Gamma \bs D_{\Sig}$ is Hausdorff.

$(2)$ When $\Gamma$ is neat, $\Gamma \bs D_{\Sig}$ is a log manifold ({\rm \cite{KNU2} III 1.1.5}). 
\end{sbthm}

Here we say that $\Gamma$ is {\it neat} if there is a faithful $V\in \Rep(G)$ such that for any $\gamma\in \Gamma$, the subgroup of $\C^\times$ generated by all eigenvalues of $\gamma: V_{\C}\to V_{\C}$ is torsion-free. 

\begin{sbpara} 
  The outline of the proof is as follows. 
  As in \cite{KNU2}, we can define various spaces $D_{\SL(2)}$, $D_{\BS}$, $E_{\sig}$ etc., and have the theory of CKS map. 
  Then, as in \cite{KNU2}, 
by using the CKS map, good properties of 
$\Gamma \bs D_{\Sig}$ are deduced from those of the space of $\SL(2)$-orbits $D_{\SL(2)}$, which  
reduce to 
the $\bR$-polarizable version of \cite{KNU2}. 
  We remark that what were shown in \cite{KNU2} by using $\bQ$-polarizations still hold 
under $\bR$-polarizations (\ref{pol2}) .

\end{sbpara}

\begin{sbthm}\label{t:main}
When $\Gamma$ is neat, 
$\Gamma \bs D_{\Sig}$ represents the functor 
$S\mapsto \{$isomorphism class of 
$G$-LMH over $S$ with a $\Gamma$-level structure of $\Sigma$-type$\}$.
\end{sbthm} 

The proof of \ref{t:main} is similar to the proof of \cite{KNU2} III, 2.6.6.

\section{Examples}
\label{s:ex}
  We discuss four examples of $D$ to which our theory can be applied so that 
we can give $D_{\Sig}$ for these $D$. 

\subsection{Usual period domains}
  We explain that the classical Griffiths domains \cite{Gr} and their mixed Hodge generalization in \cite{U} are 
essentially regarded as 
special cases of the period domains of this paper. 
  In this case, our partial compactifications essentially coincide with those in \cite{KNU2} III.

  Let $\Lambda=(H_0, W, (\langle \;,\;\rangle_w)_w, (h^{p,q})_{p,q})$ be as usual as in \cite{KNU2} III. 
  Let $G$ be the subgroup of $\Aut(H_{0,\Q}, W)$ consisting of elements which induce {\it similitudes for $\langle\;,\;\rangle_w$} for each $w$. That is, 
$
G:= \{g\in \Aut(H_{0,\Q}, W)\;|\;$ for any $w$, there is a $t_w\in {\mathbb G}_m$ such that 
$\langle gx,gy\rangle_w = t_w\langle x, y \rangle_w$ for any $x,y \in \gr^W_w\}.
$
  Let $G_1:=\Aut(H_{0,\Q}, W, (\langle \;,\;\rangle_w)_w) \subset G$.

  Let $D(\Lambda)$ be the period domain of \cite{U}. 
Then $D(\Lambda)$ is identified with an open and closed part 
 of  $D$ 
in this paper as follows. 

Assume that $D(\Lam)$ is not empty and fix an $\br\in D(\Lambda)$. Then the Hodge decomposition of  $\gr^W\br$  induces 
$h_0: S_{\C/\R} \to (G/G_u)_\R$. (We have
$\langle [z]x, [z]y\rangle_w= |z|^{2w}\langle x, y\rangle_w$ for $z\in \C^\times$ (see \ref{SCR} for $[z]$).) 
  Consider the associated period domain $D$ (\ref{D}). Then $D$ is a finite disjoint union of $G_1(\R)G_u(\C)$-orbits which are open and closed in $D$. 
  Let $\cD$ be the $G_1(\R)G_u(\C)$-orbit in $D$ consisting of points whose associated homomorphisms $S_{\C/\R}\to (G/G_u)_\R$  are $(G_1/G_u)(\R)$-conjugate to $h_0$. 
  Then the map $H\mapsto H(H_{0,\Q})$ gives a $G_1(\R)G_u(\C)$-equivariant  isomorphism $\cD\overset{\cong}\to D(\Lambda)$. 
  
\subsection{Mixed Mumford--Tate domains}
\label{ss:MT}
  Let $H$ be a MHS whose $\gr^W$ are $\R$-polarizable.

  The {\it Mumford--Tate group $G$ of $H$} is 
the Tannaka group (cf.\ \cite{Mi}) of the Tannaka category generated by $H$ 
(cf.\ \cite{A}).
  Explicitly, it is the smallest $\Q$-subgroup $G$ of $\Aut(H_\Q)$ such that $G_\R$ contains the image of the homomorphism $h:S_{\C/\R} \to \Aut(H_\R)$ and such that $\Lie(G)_\R$ contains $\delta$. 
  Here $h$ and $\delta$ are determined by the canonical splitting of $H$ (\cite{KNU2} II 1.2).
  In the case where $H$ is pure. 
$G$ is the smallest $\Q$-subgroup of $\Aut(H_\Q)$ such that $G_\R$ contains the image of $S_{\C/\R}\to \Aut(H_\R)$. 

  The {\it Mumford--Tate domain associated to $H$} is defined as the period domain $D$ associated to $G$ and $h_0: S_{\C/\R}\to (G/G_u)_\R$ 
which is defined by $\gr^WH$.

  In the pure case, our $\Gamma \bs D_{\Sig}$ is essentially the same as the one by Kerr--Pearlstein (\cite{KP}).

\subsection{Mixed Shimura varieties}
  See \cite{Mi} for the generality of mixed Shimura varieties. 
  This is the case where the universal object satisfies Griffiths transversality. $\gr^W_w\Lie(G)$ should be $0$ unless $w=0, -1, -2$. The $(p, q)$-Hodge component of $\gr^W_w\Lie(G)$ for $w=0$ (resp.\ $w=-1$, resp.\ $w=-2$) should be $0$ unless $(p,q)$ is $(1,-1)$, $(0,0)$, and $(-1,1)$ (resp.\ $(0, -1)$ and $(-1, 0)$, resp.\ $(-1,-1)$). (If this condition is satisfied by one point of $D$, it is satisfied by all points of $D$.)

In the case of PEL 
(polarizations, endomorphisms, and level
structures) type, toroidal compactifications of universal abelian varieties are expressed as $\Gamma \bs D_{\Sig}$.

\subsection{Higher Albanese manifolds}

The higher Albanese manifold (see \cite{HZ}) can be explained by using $D$ of this paper so that we can construct a toroidal partial compactification of it by the method of this paper.

\begin{sbpara}\label{hab1}
Let $X$ be a connected smooth algebraic variety over $\C$. Fix $b\in X$. 
Let $\Gamma$ be a quotient group of $\pi_1(X, b)$ and assume that $\Gamma$ is a torsion-free nilpotent group.

Let $\cG$ be the unipotent algebraic group over $\Q$ whose Lie algebra is defined as follows. Let $I$ be the augmentation ideal $\text{Ker}(\bQ[\Gamma]\to \bQ)$. Then $\Lie(\cG)$ is the $\Q$-subspace of $\varprojlim_n \Q[\Gamma]/I^n$ generated by all $\log(\gamma)$ ($\gamma \in \Gamma$). The Lie product of $\Lie(\cG)$ is defined by $[x,y]= xy-yx$. We have $\Gamma \subset \cG(\Q)$.

Then $\Lie(\cG)$ has a natural MHS  of weights $\leq -1$ such that $\gr^W_w\Lie(\cG)$ are polarizable and such that the Lie product $\Lie(\cG) \otimes \Lie(\cG) \to \Lie(\cG)\; ;\;x\otimes y \mapsto [x,y]$ is a homomorphism of MHS.

 \end{sbpara}
 \begin{sbpara}\label{hab2}
Let $\cC_{X,\Gamma}$ be the category of variations of $\Q$-MHS $\cH$ on $X$ satisfying the following conditions.

(i) For any $w\in \Z$, 
$\gr^W_w\cH$ is a constant  polarizable Hodge structure.

(ii) $\cH$ is good at infinity in the sense of \cite{HZ} (1.5). 

(iii) The monodromy action of $\pi_1(X, b)$ on $\cH_{\Q,b}$ (which is unipotent under (i)) factors through $\Gamma$.

For an object $\cH$ of $\cC_{X,\Gamma}$, we have the fiber $\cH(b)$ at $b$ which is a $\Q$-MHS. The $\Q$-vector space $\cH(b)_\Q$ has a unipotent linear action of $\Gamma$ and hence an action of the Lie algebra $\Lie(\cG)$.

The theorem of Hain--Zucker says that the functor $\cH \mapsto \cH(b)$ gives an equivalence of categories
$$\cC_{X,\Gamma} \overset{\simeq}\to \cC_{X,\Gamma}',$$
where $\cC_{X,\Gamma}'$ is the category of $\Q$-MHS $H$ with polarizable $\gr^W$ endowed with an action of the Lie algebra $\Lie(\cG)$ on $H_{\Q}$ such that  $\Lie(\cG) \otimes H \to H$ is a homomorphism of MHS.

\end{sbpara}

\begin{sbpara}\label{hab3} The higher Albanese manifold $A_{X,\Gamma}$ of $X$ for $\Gamma$ is as follows. 
Let $F^0\cG(\C)$ be the algebraic subgroup of $\cG(\C)$ over $\C$ corresponding to the Lie subalgebra $F^0\Lie(\cG)_\C$ of $\Lie(\cG)_\C$. Define
$$A_{X,\Gamma}:= \Gamma\bs \cG(\C)/F^0\cG(\C).$$

In the case where $\Gamma$ is $H_1(X, \Z)/(\text{torsion})$ regarded as a quotient group of $\pi_1(X,b)$, $A_{X,\Gamma}$ coincides with the Albanese variety 
$\Gamma\bs H_1(X,\C)/F^0H_1(X,\C)$ of $X$. 

We will give an understanding of $A_{X,\Gamma}$ by using $D$ of this paper (\ref{und8}). 

We will describe the functor represented by $A_{X,\Gamma}$ (Theorem \ref{thm3} (1)). 

\end{sbpara}

\begin{sbpara}\label{und2} Let $Q$ be the Mumford--Tate group (\ref{ss:MT}) of the MHS $\Lie(\cG)$ (\ref{hab1}). The action of $Q$ on $\Lie(\cG)$ induces an action of $Q$ on $\cG$. 
By using this action, define the semidirect product $G$ of $Q$ and $\cG$ with an exact sequence $1\to \cG\to G\to Q\to 1$. We have $\cG\subset G_u$. We have $h_0: S_{\C/\R}\to (Q/Q_u)_\R= (G/G_u)_\R$ given by the Hodge decomposition of $\gr^W\Lie(\cG)$.

Then $(G, \Gamma)$ satisfies the condition in \ref{Gamma}, and $\Gamma$ is a neat subgroup of $G(\Q)$. 

Let $D_G$ (resp.\ $D_Q$) be the period domain $D$ for $G$ (resp.\ $Q$) and $h_0$. 
We have a canonical map $\Gamma \bs D_G\to D_Q$ induced by the canonical homomorphism $G\to Q$. 

\end{sbpara}

\begin{sbpara}\label{und7}

The equivalence in \ref{hab2} gives an inclusion functor 
$$\Rep(G)\overset{\subset}\to  \cC_{X,\Gamma}.$$ 

\end{sbpara}

\begin{sbpara}\label{und8}

Let $b_G\in D_G$ be the element induced by $\cC_{X, \Gamma} \to \text{MHS}\;;\; \cH \mapsto \cH(b)$ and \ref{und7}. 
Let $b_Q\in D_Q$ be the image of $b_G$.

 The map $G_u(\C) \to D_G\; ;\; g \mapsto gb_G$ induces an isomorphism from $A_{X,\Gamma}=\Gamma\bs \cG(\C)/F^0\cG(\C)$ to the inverse image of $b_Q$ in $\Gamma\bs D_G\to D_Q$. 
 
\end{sbpara}

\begin{sbpara} Let $\Sig$ be a weak fan in $\Lie(G)$ such that $\sig\subset \Lie(\cG)_\R$ for any $\sig\in \Sig$. We have a canonical morphism $\Gamma \bs D_{G,\Sig}\to D_Q$ induced by the homomorphism $G\to Q$.

Define the toroidal partial compactification $
{A}_{X, \Gamma,\Sig}$ of $A_{X,\Gamma}$ as the subspace of $\Gamma \bs D_{G,\Sig}$ defined to be the inverse image of $b_Q$. We can endow $
{A}_{X, \Gamma, \Sig}$ with a structure of a log manifold such that for any object $S$ of $\cB(\log)$,  $\text{Mor}(S, 
{A}_{X,\Gamma,\Sig})$ coincides with the set of all morphisms $S \to \Gamma \bs D_{G,\Sig}$ whose images in $D_Q$ are $b_Q$.

\end{sbpara}

\begin{sbpara}\label{moduli} Define contravariant functors $$\cF_{\Gamma}, 
{\cF}_{\Gamma,\Sig}:\cB(\log) \to(\text{Set})$$
as follows. 

$
{\cF}_{\Gamma,\Sig}(S)$ is the set of isomorphism classes of pairs $(H,\lambda)$, where $H$ is an exact $\otimes$-functor
$\cC_{X,\Gamma} \to \LMH(S)$ and  $\lambda$ is a global section of the sheaf $\Gamma\bs {\cal I}$ on $S^{\log}$, where $\cal I$ is the sheaf of functorial $\otimes$-isomorphisms $H(\cH)_{\Q} \overset{\cong}\to \cH(b)_{\Q}$ of $\Q$-local systems, satisfying the following conditions (i) and (ii).

 (i) For any $\Q$-MHS $h$, we have a functorial $\otimes$-isomorphism $H(h_X) \cong h$ such that the induced isomorphism of local systems $H(h_X)_\Q\cong h_\Q=h_X(b)_\Q$ belongs to $\lambda$. Here $h_X$ denotes the constant variation of $\Q$-MHS over $X$ associated to $h$.

 (ii) The following (ii-1) and (ii-2) are satisfied for any $s\in S$ and any $t\in s^{\log}$. Let $\tilde \lambda_t: H(\cH)_{\Q,t}\cong \cH(b)_\Q$ be a  functorial $\otimes$-isomorphism which belongs to $\lambda_t$.

  (ii-1) There is a $\sig\in \Sig$ such that the logarithm of the action of $\Hom((M_S/\cO^\times_S)_s, \N)\subset \pi_1(s^{\log})$ on $H_{\Q,t}$ is contained, via $\tilde \lambda_t$, in $\sigma \subset \Lie(\cG)_\R$.

  (ii-2) Let $\sig\in \Sig$ be the smallest cone which satisfies (ii-1) and let $a: \cO_{S,t}^{\log}\to \C$ be a ring homomorphism which induces the evaluation $\cO_{S,s}\to \C$ at $s$. Then for each $\cH\in \cC_{X,\Gamma}$, $(\sig, \tilde \lambda_t(a(H(\cH))))$ generates a nilpotent orbit in the sense of \cite{KNU2} III 2.2.2.

 We define a subfunctor $\cF_{\Gamma}\subset 
{\cF}_{\Gamma,\Sig}$ by replacing $\LMH(S)$ in the above definition of $
{\cF}_{\Gamma,\Sig}$ by the full subcategory $\text{MHS}(S)\subset \LMH(S)$ consisting of objects with no degeneration. That is, $\text{MHS}(S)$ is the category of analytic families of mixed Hodge structures parametrized by $S$. (In this case, the above condition (ii) is empty. Objects of $\text{MHS}(S)$ need not satisfy Griffiths transversality.) 
 
 \end{sbpara}

\begin{sbthm}\label{thm3}

$(1)$ The functor $\cF_{\Gamma}$ is represented by $A_{X,\Gamma}$. 

$(2)$ The functor $
{\cF}_{\Gamma,\Sig}$ is represented by $
{A}_{X, \Gamma,\Sig}$.
\end{sbthm}

This theorem is proved in the following way. 
By 
\ref{und7}, \ref{und8}, and Theorem \ref{t:main}, we have a map from $\cF_{\Gamma}(S)$ (resp.\ $
{\cF}_{\Gamma,\Sig}(S)$) to the fiber $\text{Mor}(S, A_{X,\Gamma})$ (resp.\ $\text{Mor}(S, 
{A}_{X,\Gamma,\Sig})$) of $\text{Mor}(S, \Gamma\bs D_G) \to \text{Mor}(S,  D_Q)$ (resp.\ $\text{Mor}(S, \Gamma\bs D_{G,\Sig}) \to \text{Mor}(S,  D_Q)$)
over $b_Q$. We can show that this map is a bijection.

\begin{sbpara} The {\it higher Albanese map} $X\to A_{X,\Gamma}$ 
 corresponds in \ref{thm3} (1) to 
 the evident functor $\cC_{X,\Gamma} \to \text{MHS}(X)$. 
 When this 
extends to a morphism $\overline{X} \to 
{A}_{X,\Gamma,\Sig}$ for some complex analytic manifold $\overline{X}$ which contains $X$ as a dense open subset such that the complement $\overline{X}\smallsetminus X$ is a divisor with normal crossings, 
this extended higher Albanese map corresponds in \ref{thm3} (2) to the 
  inclusion functor $\cC_{X,\Gamma} \to \LMH(\overline{X})$. 
\end{sbpara}

{\bf Acknowledgments.} 
K.\ Kato was 
partially supported by NFS grants DMS 1001729 and DMS 1303421.
C.\ Nakayama was 
partially supported by JSPS Grants-in-Aid for Scientific Research (C) 22540011 and (B) 23340008.
S.\ Usui was 
partially supported by JSPS Grants-in-Aid for Scientific Research (B) 23340008.

\end{document}